\documentclass[11pt]{article}
\usepackage{amssymb,latexsym,amsmath}
\usepackage{setspace}
\long\def\symbolfootnote[#1]#2{\begingroup%
\def\thefootnote{\fnsymbol{footnote}}\footnote[#1]{#2}\endgroup}
\input{amssym.def}
\input{amssym}
\newtheorem{definition}{Definition}
\newtheorem{theorem}{Theorem}
\newtheorem{remark}{Remark}

\newtheorem{example}{Example}
\newtheorem{lemma}{Lemma}

\renewcommand{\baselinestretch}{1.25}
\begin{document}
%\doublespace
%{\noindent  Draft of  \today:}
\begin{center}
 \Large{\sc  Intrinsic   posterior regret gamma-minimax estimation for the exponential family of distributions }
 %\footnote{draft of \today}
\end{center}

\begin{center}
{Mohammad Jafari Jozani$^{a,}$\footnote{Corresponding author: m$_{-}$jafari$_{-}$jozani@umanitoba.ca} and  Nahid Jafari Tabrizi$^b$}

\vspace{0.3cm}

{\it $^a$ Department of Statistics,  University of Manitoba, 
Winnipeg, MB, CANADA, R3T 2N2.

{\it $^b$ Department of Mathematics,  Islamic Azad University-Karaj Branch,  Karaj, IRAN.}}
\end{center}

\noindent {\bf \sc ABSTRACT: } 
 \noindent \small
In practice, it is desired to have estimates that are invariant  under  reparameterization. The invariance property of the estimators helps  to formulate a unified solution to the underlying estimation  problem. In robust Bayesian analysis, a frequent criticism is that the optimal estimators are not invariant under smooth reparameterizations. This paper considers the problem of posterior regret gamma-minimax (PRGM) estimation of the natural parameter of the exponential family of distributions under intrinsic loss functions. We show that under  the class of Jeffrey's Conjugate Prior (JCP) distributions, PRGM estimators are invariant to smooth one-to-one reparameterizations. We apply our results to several distributions and different  classes of JCP, as well as  the usual  conjugate  prior distributions. We observe  that, in many cases,  invariant PRGM estimators in the class of JCP distributions can be obtained by some modifications of PRGM estimators in  the usual class of conjugate priors. 
 Moreover, when the class of priors are convex or  dependant on a hyper-parameter belonging to  a connected set, we show that the PRGM estimator under the intrinsic loss function could be Bayes with respect to a prior distribution in the original prior class. Theoretical results  are supplemented  with several examples and illustrations.


%\noindent MSC:  62C10, 62F10
\normalsize
\section{\sc Introduction}

%This papers considers robust Bayesian analysis under intrinsic loss functions.
Suppose $x$ is a realization of a  random sample $X$  with a sampling model given by a family of densities $\{ f(\cdot| \theta):  \theta\in\Theta\}$ with respect to a $\sigma$-finite measure $\nu$ on a sample space $\chi$ where $\theta$ is the unknown parameter of interest  with $\theta\in\Theta$. Let $\pi(\cdot)$ be a prior distribution on $\Theta$ and $ \pi(\cdot|x)$ denote the posterior distribution of $\theta$ given  $x$.  
In standard Bayesian analysis,  one needs to specify the true prior distribution $\pi(\cdot)$. However, in practice,    elicitation of the true prior distribution  can never be done without error. Hence,  we  usually need to consider   a  class $\Gamma$ of prior distributions which reflect (approximately) true prior beliefs, i.e., the true prior distribution $\pi(\cdot)$ is an unknown element of $\Gamma$. 
 Robust Bayesian analysis is designed   to acknowledge such
a prior uncertainty by considering the class $\Gamma$ of plausible prior distributions   instead of a single prior distribution $\pi$  and studying
the corresponding range of Bayesian solutions.  See    Berger (1994) and Rios Insua and Ruggeri (2000) for more details. One may also attempt to determine an optimal estimator 
$\delta$ by  minimizing some measures of robustness. Several criteria  have been proposed for the selection of procedures in robust Bayesian studies. In this paper,  we study the  maximal posterior regret method (e.g., Rios Insua and Ruggeri, 2000;  Rios Insua et al., 1995) to obtain the posterior regret
gamma-minimax  (PRGM) estimator of  the unknown  parameter  for  the one-parameter exponential family of distributions.  The PRGM criterion has been used recently by many people from both theoretical and practical points of view. For example,  G\'omez-D\'eniz (2009) investigated the use of PRGM for credibility premium estimation in Actuarial Science,  Boraty\'nska (2002, 2006) in insurance for collective risk model  analysis, and   Jafari Jozani and Parsian (2008)  in  statistical inference  based on record data.  

For an observed value $x$, a prior distribution $\pi$ and the corresponding posterior distribution $\pi(\cdot|x)$, we denote the  posterior risk of  an estimate $\delta(x)$ of the unknown parameter $\theta$  under  $L(\theta, \delta)$ by  $r(x, \delta)=E[L(\theta, \delta(x))| x]$. The Bayes estimator of $\theta$ under the loss  function $L(\theta, \delta)$ is then given by a $\delta_{\pi}(X)$ such that  $r(x, \delta_{\pi})=\inf_{\delta} r(x, \delta) $.

\begin{definition}
 The {\rm PRGM} estimator of $\theta$ under the loss function $L(\theta, \delta)$ and a class $\Gamma$ of prior distributions is defined as an estimator  $\delta_{PR}$  such that 
\begin{align}\label{PRGM}
\sup_{\pi\in\Gamma}\rho(\delta_{\pi}(x),
\delta_{PR}(x))=\inf_{\delta}\sup_{\pi\in\Gamma}\rho(\delta_\pi(x), \delta(x)),
\end{align}
where $\rho(\delta_{\pi}, \delta)=r(x, \delta)-r(x, \delta_{\pi})$ is the
posterior regret measuring the loss entailed in choosing the action $\delta(x)$ instead of the
optimal Bayes action $\delta_{\pi}(x)$ (under prior $\pi$ and  loss $L$).
\end{definition}

 In this paper,  we study  the construction of  PRGM estimators  under the so-called \textit{intrinsic loss} functions.    These loss functions shift attention from the distance between the estimator $\delta$ and the true parameter value $\theta$, to the more relevant distance  between statistical models they label. More specifically, the intrinsic loss of using $\delta$ as a proxy for $\theta$ is the intrinsic distance between the true model $f(x| \theta)$ and the model $f(x| \delta)$ when $\theta=\delta$, that is 
\begin{align} \label{int-def}L(\theta, \delta)= d\left( f(x| \theta) , f(x| \delta)\right),\end{align}
%$$ \label{int-def}L(\theta, \delta)= d\left( f(x| \theta) , f(x| \delta)\right)$$
where $d(\cdot, \cdot)$ is a suitable distance measure. 
In practice, intrinsic loss functions could be used as benchmark losses when the utility function related to the underlying statistical problem cannot be obtained by practitioners. A  desired property of intrinsic loss functions is that they are invariant under one-to-one smooth  reparameterizations. The invariance property of intrinsic loss functions  provides a very convenient tool for statistical application. We show that,  under suitable conditions, intrinsic loss functions could be used to formulate a unified set of solutions to the problem of  PRGM estimation of the unknown parameter of the exponential family of distributions which is consistent under reparameterization, a rather obvious requirement, which unfortunately many statistical methods fail to satisfy.  

In Section 2,  we obtain   the PRGM estimator of the natural parameter $\theta$ of the exponential family of distributions  under the intrinsic loss function \eqref{int-def} when $d(\cdot, \cdot)$  is chosen to be the Kullback-Leibler distance.   We consider  different classes of conjugate priors on the natural parameter $\theta$ and show how to obtain the PRGM estimator of $\theta$ in each class.  The results are very general and provide  an automated and unified  solution to the PRGM estimation of the unknown parameter of  the exponential family of distributions under different loss functions,  including, but not limited to,    quadratic, LINEX, entropy and  Stein loss functions. 

In Bayesian statistical analysis, as pointed out by Gelman (2004),  transformations of the  parameter typically suggest  new families of prior distributions.  Therefore, the usual robust Bayesian inferences  are not invariant under reparameterizations. For example, if  $\delta_{PR}(X)$ is the PRGM estimator of $\theta$, then it is not necessarily true that $h(\delta_{PR}(X))$ is the PRGM estimator of $\eta=h(\theta)$,  when $h$ is a one-to-one smooth function. A solution to this problem is proposed in   Section 3. To this end, we obtain invariant PRGM estimators of $\theta$ under the intrinsic loss function and   different classes of Jeffrey's Conjugate Prior (JCP) distributions.  We show that the resulting PRGM estimates are invariant under one-to-one smooth transformations of  $\theta$. 
 Theoretical results are augmented  with several examples and
illustrations.
In Section 4,  we provide some general results showing that,   under general conditions,  PRGM and intrinsic PRGM estimators are Bayes with respect to  prior distributions in the underlying class of priors. We study two cases of convex classes of prior distributions as well as the case where the underlying class of priors depends on a hyper-parameter belonging to  a connected set. 
 We provide a sufficient condition under which  the PRGM  and  intrinsic PRGM estimators  are  Bayes with respect  to data independent prior distributions within the underlying class of priors.  
 Finally,  in Section 5,  we give some concluding remarks.

\section{\sc   PRGM estimation under intrinsic loss functions}\label{IPRGM-se}

Suppose $X$ is a random variable, where its distribution belongs to the one-parameter exponential family of distributions  $\mathcal{F}=\{f(x| \theta): x\in\chi \subseteq\mathbb{R}, \theta\in\Theta\subseteq\mathbb{R} \}$, with probability density function (pdf)
\begin{align}\label{exp-family}
f(x| \theta) = \beta(\theta) t(x) e^{-\theta r(x)}, 
\end{align} 
where $r(x)>0$, $\beta(\theta) t(x)>0$ and $\theta$ is the unknown real-valued natural parameter  of the model. 
 The density is considered with respect to the Lebesgue measure for continuous  and the counting measure for discrete distributions. 
 Suppose $\delta $ is an estimate of $\theta$ with both $\theta, \delta\in\Theta$. We define   the intrinsic loss function 
\eqref{int-def},   using   the Kullback-Leibler measure  between $f(x| \theta)$ and $f(x| \delta)$,  as follows
 \begin{align}\label{KL}
 L(\theta, \delta)= E_{\theta}\left[ \log\left(\frac{f(X| \theta)}{f(X| \delta)}\right)\right]= \int_{\chi} \log\left(\frac{f(x| \theta)}{f(x| \delta)}\right) f(x| \theta) d\nu(x).
 \end{align}  
Loss function \eqref{KL} can be interpreted as the expected log-likelihood ratio in favour of the true model. Thus,  the  intrinsic loss function \eqref{KL}  not only has the desired invariance property but it is also related to the relevant measure of evidence in the Neyman-Pearson Lemma.
Note that  the intrinsic loss function \eqref{KL}  is  invariant under reparameterization since  the parameters affect the loss function  only via the probability distributions they label, which are independent of the particular parameterization.   For a general reference on intrinsic losses and
additional details we refer to Robert (1996) and Bernardo (2011).

  First, we give a lemma which identifies the intrinsic loss function for the exponential family of distributions.  

 \begin{lemma} For the  exponential family   of distributions  \eqref{exp-family},  the intrinsic loss function \eqref{KL} reduces to 
 \begin{align}\label{KL-loss}
 L(\theta, \delta)= \log\left( \frac{\beta(\theta)}{\beta(\delta)}\right) + (\delta-\theta) \frac{\beta'(\theta)}{\beta(\theta)},
 \end{align}
 where $\beta'(\theta)=\frac{d}{d\theta} \beta(\theta)$.
 \end{lemma}

 %%%%%%%%%%%%%%%%%

\noindent
 Let $H(t):= {\beta'(t)}/{\beta(t)}$.  A straightforward calculation shows that the posterior risk associated with  $\delta$,  under the loss function \eqref{KL-loss}, is  
\begin{align}\label{P-rs}
r(x,\delta)=E(\log\beta(\theta) | x)-\log\beta(\delta(x))+\delta(x)
E\left(H(\theta) \big|x\right)-E\left(\theta H(\theta)\big| x\right).
\end{align}
The Bayes  estimator of $\theta$ can therefore be obtained by minimizing \eqref{P-rs} in $\delta$ as follows   
\begin{align}\label{Bayes-general}
\delta_{\pi}(X)=H^{-1}\{E\left(H(\theta)\big|X\right)\}.
\end{align}
  Following the decreasing monotone likelihood ratio property of the densities  $f(x|\theta)$ in \eqref{exp-family} in $r(X)$,  and since $E[r(X)]= H(\theta)$,  $H(\cdot)$ is a decreasing function. Therefore,    the Bayes estimator $\delta_{\pi}(X)$   is unique.
Furthermore, the posterior regret for estimating $\theta$ using $\delta$ instead of  the optimal estimator $\delta_{\pi}$  is  obtained by 
\begin{align}\label{P-rg}
\rho(\delta_{\pi},\delta)= \log\frac{\beta(\delta_{\pi})}{\beta(\delta)}+(\delta-\delta_{\pi})H(\delta_{\pi}).
\end{align}
Note that $\rho(\delta_{\pi}, \delta)$, as a function of $\delta_{\pi}$, 
decreases then increases with a   unique minimum at $\delta_{\pi}=\delta$. 
%In other word
%\begin{align}\label{derivative}
%\frac{\partial\rho(\delta_{\pi},\delta)}{\partial\delta_{\pi}}=\left\{\begin{array}{ll}
%<0,&\hbox{~~}\delta_{\pi}<\delta, \cr
%=0,&~~\delta_{\pi}=\delta, \cr
%>0,&~~\delta_{\pi}>\delta.
%\end{array}\right.
%\end{align}
%Note that the exponential family (3) has MLR property in $-r(x)$ so
%$E(-r(X))=-\frac{\beta'(\theta)}{\beta(\theta)}$ is an increasing
%function of $\theta$.
The main result of this section is given in the following theorem which  obtains  the PRGM estimator of  $\theta$   under the intrinsic loss function \eqref{KL-loss}.

\begin{theorem}\label{the1} Let $\underline{\delta}(x)=\inf_{\pi\in\Gamma}\delta_{\pi}(x)$ and
$\overline{\delta}(x)=\sup_{\pi\in\Gamma}\delta_{\pi}(x)$ and suppose that $\underline{\delta}(x)$ and $\overline{\delta}(x)$ are finite almost everywhere. The PRGM estimator of $\theta$ in the exponential family \eqref{exp-family} under
the loss function \eqref{KL-loss} and  in the class of prior distributions $\Gamma$  is given by
\begin{align}\label{PRGM-est}
\delta_{PR}(X)=\frac{\overline{\delta}(X)\, H(\overline{\delta}(X))-\underline{\delta}(X)\, H(\underline{\delta}(X))-\log\frac{\beta(\overline{\delta}(X))}{\beta(\underline{\delta}(X))}}{ H(\overline{\delta}(X))-H(\underline{\delta}(X))}.
\end{align}
\end{theorem}
\textbf{Proof:} First,  note that $$\inf_{\delta}\sup_{\pi\in\Gamma} \rho(\delta_{\pi}, \delta)= \min\left\{\inf_{\delta\leq\underline{\delta}}~\sup_{\pi\in\Gamma}~\rho(\delta_{\pi},\delta), \inf_{\underline{\delta}<\delta<\overline{\delta}}~\sup_{\pi\in\Gamma}~\rho(\delta_{\pi},\delta),~  \inf_{\delta\geq\overline{\delta}}~\sup_{\pi\in\Gamma}~\rho(\delta_{\pi},\delta) \right\}.$$
So,  we consider the following three cases:

 \noindent\textbf{Case 1. }  When   $\delta\leq\underline{\delta}$,  we have  
 $\sup_{\pi\in\Gamma}\rho(\delta_{\pi},\delta)=\rho(\overline{\delta},\delta)$.
 %for $\delta\leq\underline{\delta}<\overline{\delta}$. 
 Let
 $f_{1}(\delta)=\rho(\overline{\delta},\delta)=\log\frac{\beta(\overline{\delta})}{\beta(\delta)}+(\delta-\overline{\delta})H(\overline{\delta})$ with 
 $f'_{1}(\delta)=H(\overline{\delta})-H(\delta)<0$,
 following the decreasing  property  of $H(\cdot)$.  Hence,    $f_{1}(\delta)$
 is  a decreasing function of $\delta$ for  $\delta\leq\underline{\delta}$ and 
 $\inf_{\delta\leq\underline{\delta}}f_{1}(\delta)=f_{1}(\underline{\delta})$.
 Therefore,
 $$\inf_{\delta\leq\underline{\delta}}~\sup_{\pi\in\Gamma}~\rho(\delta_{\pi},\delta)=\rho(\overline{\delta},\underline{\delta}).$$

\noindent\textbf{Case 2.} For $\delta\geq\overline{\delta}$,  we have  
 $\sup_{\pi\in\Gamma}\rho(\delta_{\pi},\delta)=\rho(\underline{\delta},\delta)$.
% for $\underline{\delta}<\overline{\delta}\leq\delta$.
  Let
 $f_{2}(\delta)=\rho(\underline{\delta},\delta)=\log\frac{\beta(\underline{\delta})}{\beta(\delta)}+(\delta-\underline{\delta})H(\underline{\delta})$ with 
 $f'_{2}(\delta)=H(\underline{\delta})- H({\delta})>0$.  Hence, $f_{2}(\delta)$ is an increasing function of $\delta$ for  $\delta\geq\overline{\delta}$
 and 
 $\inf_{\delta\geq\overline{\delta}}f_{2}(\delta)=f_{2}(\overline{\delta})$.
 Therefore,
 $$\inf_{\delta\geq\overline{\delta}}~\sup_{\pi\in\Gamma}~\rho(\delta_{\pi},\delta)=\rho(\underline{\delta},\overline{\delta}).$$

\noindent\textbf{Case 3.} If
$\underline{\delta}<\delta<\overline{\delta}$,  then
$\sup_{\pi\in\Gamma}\rho(\delta_{\pi},\delta)=\max\{\rho(\overline{\delta},\delta),\rho(\underline{\delta},\delta)\}$.
Let
$f_{3}(\delta)=f_{1}(\delta)-f_{2}(\delta)$ where    $f'_{3}(\delta)= H(\overline{\delta})- H(\underline{\delta})<0$.
 Since $f_{3}(\delta)$ is a decreasing function of $\delta$
 with $f_{3}(\overline{\delta})<0$ and $f_{3}(\underline{\delta})>0$, 
%  $f_{3}(\overline{\delta})=-(\log\frac{B(\underline{\delta})}{B(\delta)}+(\delta-\underline{\delta})H(\underline{\delta})<0$ and $f_{3}(\underline{\delta})=\log\frac{B(\overline{\delta})}{B(\delta)}+(\delta-\overline{\delta})
% H(\overline{\delta})>0$,
  there exists a unique
 $\delta^{*}\in(\underline{\delta},\overline{\delta})$  (as  the root of
 $f_{3}(\delta)=0$) such that
 $\rho(\underline{\delta},\delta^{*})=\rho(\overline{\delta}, \delta^{*})$.
 Hence, for $\underline{\delta}<\delta<\delta^{*}$,
 $\sup_{\pi\in\Gamma}\rho(\delta_{\pi},\delta)=\rho(\overline{\delta},\delta)$
 and for $\delta^{*}<\delta<\overline{\delta}$, $\sup_{\pi\in\Gamma}\rho(\delta_{\pi},\delta)=
 \rho(\underline{\delta},\delta)$.
Note that, for $\underline{\delta}<\delta<\overline{\delta}$,
$\rho(\overline{\delta},\delta)$ is  a decreasing function in
$\delta$ with 
$\inf_{\underline{\delta}<\delta<\delta^{*}}\sup_{\pi\in\Gamma}~\rho(\delta_{\pi},\delta)=\rho(\overline{\delta},\delta^{*})$
and $\rho(\underline{\delta},\delta)$ is an increasing function in
$\delta$ with 
$\inf_{\delta^{*}<\delta<\overline{\delta}}\sup_{\pi\in\Gamma}~\rho(\delta_{\pi},\delta)=\rho(\underline{\delta}, \delta^{*})$.
Therefore, 
$$\inf_{\underline{\delta}<\delta<\overline{\delta}}\sup_{\pi\in\Gamma}~\rho(\delta_{\pi},\delta)=\rho(\underline{\delta}, \delta^*)
=\rho(\overline{\delta},\delta^{*}).$$
\noindent Following  the  above cases, we conclude that 
$$\inf_{\delta\in D}\sup_{\pi\in\Gamma}~\rho(\delta_{\pi},\delta)=\inf_{\underline{\delta}<\delta<\overline{\delta}}~
\sup_{\pi\in\Gamma}~\rho(\delta_{\pi},\delta)=\rho(\underline{\delta}, \delta^*)
=\rho(\overline{\delta},\delta^{*}). $$ That is, the PRGM estimator of $\theta$  is given by 
$\delta_{PR}=\delta^{*}\in (\underline\delta, \overline\delta)$, as   the solution of 
$$\log\frac{\beta(\overline{\delta})}{\beta(\underline{\delta})}+\delta_{PR}\big(H(\overline{\delta})-H(\underline{\delta})\big)+ \underline{\delta}H(\underline{\delta})-\overline{\delta}H(\overline{\delta})=0,$$
 in $\delta_{PR}$ which results in the estimator  \eqref{PRGM-est}. $\Box$\\
%Therefore
%$$\delta^{PRGM}=\frac{\overline{\delta}\frac{B'(\overline{\delta})}
% {B(\overline{\delta})}-\underline{\delta}\frac{B'(\underline{\delta})}{B(\underline{\delta})}-\log\frac{B(\overline{\delta})}
% {B(\underline{\delta})}}{\frac{B'(\overline{\delta})}{B(\overline{\delta})}-\frac{B'(\underline{\delta})}{B(\underline{\delta})}}
% $$.

%%%%%%%%%%%%%%%%%%%%%%%

We give  some  applications of Theorem \ref{the1}.

 \begin{example}\label{normal} (Normal distribution). Suppose $X\sim N(\mu,1)$ is a normally distributed  random variable  with unknown parameter $\mu\in\mathbb{R}$ and pdf
$f(x|\mu)=\frac{1}{\sqrt{2\pi}}e^{-\frac{1}{2}(x-\mu)^{2}}$, $-\infty<x<\infty.$
The pdf $f(x |\mu)$ belongs to the exponential family \eqref{exp-family} with $\theta=\mu$, 
and $\beta(\theta)=e^{-\frac{\theta^{2}}{2}}$.  Also,  $H(\theta)=-\theta$, and   the
intrinsic  loss function \eqref{KL-loss} reduces to $L(\theta, \delta)
=\frac{1}{2}(\delta-\theta)^{2}$ which is essentially the  usual squared error loss function.  Let   $\underline\delta$ and $\overline\delta$  be defined as  in Theorem \ref{the1}. Using \eqref{PRGM-est}, subject to the existence of $\overline{\delta}$ and $\underline{\delta}$,  the PRGM estimator of $\theta$ in  the class $\Gamma$ of prior distributions is given by 
$$\delta_{PR}(X)=\frac{1}{2}(\overline{\delta}(X)+\underline{\delta}(X)),$$
  which is also obtained in  Rios Insua et al.\ (1995) as well as Berger (1994).
\end{example}

\begin{example}\label{exponential-dist} (Exponential distribution).   Suppose $X\sim Exp(\sigma)$ is an exponential random variable with  pdf $f(x| \sigma) =\frac{1}{\sigma} e^{- x/\sigma}$, $x>0$, where $\sigma>0$  is the unknown parameter. The pdf $f(x| \sigma)$  belongs to the 
exponential family \eqref{exp-family} with  $\theta= \frac{1}{\sigma}$,  and   $\beta(\theta)=\theta$. In this case, $H(\theta)={\theta}^{-1}$,  and  the
intrinsic  loss function \eqref{KL-loss}  reduces to the Stein  loss   $$L(\theta, \delta)=
 \frac{\delta}{\theta}-\log
\frac{\delta}{\theta}-1.$$ Using  \eqref{PRGM-est}, subject to the existence of $\overline{\delta}$ and $\underline{\delta}$,
the PRGM estimator of $\theta$ under the Stein  loss function is given by 
$$\delta_{PR}(X)=\frac{\log\frac{1}{\overline{\delta}(X)} - \log\frac{1}{\underline{\delta}(X)}}{\frac{1}{\overline{\delta}(X)}-\frac{1}{~\underline{\delta}(X)}}.$$
The PRGM estimator of  $\sigma$ is also obtained in Example \ref{ex-exp-bayes}. 
 \end{example}

\begin{example}\label{binomial-dist} (Binomial distribution).  Suppose  $X\sim Bin(n, p)$ is a binomial random variable  with probability mass function (pmf)  $f(x|p )=\binom{n}{x} p^{x}(1-p)^{n-x}$, where $n$ is known, 
$x=0, 1, \ldots, n,$ and $p\in[0, 1]$ is the  unknown parameter.   The pmf  $f(x|p)$ is a member of the exponential family \eqref{exp-family} with  $\theta=\log(\frac{1-p}{p})$ and $\beta(\theta) = (1+ e^{-\theta}) ^{-n}$. We also have  $H(\theta)= \frac{n}{1+ e^{\theta}}$ which results in the    intrinsic loss function \begin{align}\label{binom-loss}L(\theta, \delta) =n\left\{  \log\left(\frac{e^{\theta}}{e^{\delta}}\cdot\frac{1+ e^{\delta}}{1+ e^{\theta}} \right) +  \frac{\delta-\theta}{1+ e^{\theta}} \right\}. \end{align} 
Using  \eqref{PRGM-est}, subject to the existence of $\overline{\delta}$ and $\underline{\delta}$,  the PRGM estimator of $\theta$  is given by 
\begin{align}\label{prgm-binom}\delta_{PR}(X) = \frac{\frac{\overline{\delta}(X) }{1+ e^{ \overline{\delta}(X)} }- \frac{\underline{\delta}(X)}{1+ e^{ \underline{\delta}(X)} }  -\log\left\{ \frac{e^{ \overline{\delta}(X)}}{e^{ \underline{\delta}(X)}} \frac{1+ e^{ \underline{\delta}(X)}}{1+ e^{ \overline{\delta}(X)} }\right\} } {\frac{1}{1+ e^{ \overline{\delta}(X)}}- \frac{1}{1+ e^{ \underline{\delta}(X)} }}.\end{align}
In Example  \ref{binom-prgm-invariant}, we obtain the  PRGM estimator of $p$. \end{example}

 %%%%%%%%%%%%%%%%%%%
 We now consider the PRGM estimation of  $\theta$ under conjugate classes of prior distributions. 
 For the exponential family \eqref{exp-family}  and   a conjugate  prior distribution    
 \begin{align}\label{prior-density}
 \pi_{\alpha, \lambda}(\theta) \propto \{ \beta(\theta)\} ^{\alpha}\,  e^{-\theta\, \lambda}, 
 \end{align} the posterior distribution  is given by  
 $\pi(\theta| x) \propto \{ \beta(\theta)\}^ {1+ \alpha} \, e^{-(\lambda+ r(x)) \theta}$, and  $\pi(\theta| x) =\pi_{\alpha+1, \lambda+ r(x)}(\theta)$.  Also,  as established by Diaconis and Ylvisaker (1979),  
 $E[H(\theta)|x]= \frac{\lambda+ r(x)}{\alpha+1}$. Now,    the Bayes estimator of $\theta$ under the intrinsic loss function \eqref{KL-loss} is obtained by (e.g., Bernardo and Smith (1994), Robert (1996) and Gutierrez-Pena(1992))
 \begin{align}\label{Bayes}
 \delta_{\pi}(X)= H^{-1} \left( \frac{\lambda+ r(X)}{\alpha+1}\right).
 \end{align}
  Furthermore, the posterior regret  for   estimating $\theta$ with  $\delta(x)$ is  $\rho(\delta_{\pi}, \delta) = \log\frac{\beta(\delta_{\pi}(x))}{\beta(\delta(x))}+ (\delta(x)-\delta_{\pi}(x)) \frac{\lambda+ r(x)}{\alpha+1}$.
 Now, suppose that the prior distribution $\pi_{\alpha, \lambda}$  belongs to the following class of conjugate prior   distributions: 
\begin{align*}
\Gamma = \{ \pi_{\alpha, \lambda}(\theta) ~: ~ \alpha\in [\alpha_1, \alpha_2] , ~ \lambda\in [\lambda_1, \lambda_2] \}, 
\end{align*}  
 with suitable choices of $\alpha_1<\alpha_2$ and $\lambda_1< \lambda_2$ leading to proper posterior distributions for $\theta$. A straightforward calculation    shows that $ H(\underline{\delta}(x)) = \frac{\lambda_2+ r(x)}{\alpha_1+ 1}$ and   $H(\bar{\delta}(x))= \frac{\lambda_1+ r(x)}{\alpha_2 + 1}$.
  Hence, we can state the following result.

 \begin{lemma} Suppose $U(t) = H^{-1}(t)$ with $H(t) = \beta'(t) / \beta(t)$. The PRGM estimate of $\theta$ for the exponential family \eqref{exp-family} under the intrinsic loss function \eqref{KL-loss} and in the class  $\Gamma$  of prior distributions is given by  
 \begin{align}\label{gamma3}
  \delta^{\Gamma}_{PR}(x)= \frac{ \frac{\lambda_1+ r(x)}{\alpha_2 + 1}U\left( \frac{\lambda_1+ r(x)}{\alpha_2 + 1} \right)  - \frac{\lambda_2+ r(x)}{\alpha_1 + 1} U\left( \frac{\lambda_2+ r(x)}{\alpha_1 + 1}\right)-\log \left(\frac{\beta(U( \frac{\lambda_1+ r(x)}{\alpha_2 + 1}) )}{\beta(U( \frac{\lambda_2+ r(x)}{\alpha_1 + 1})) } \right)  }{  \frac{\lambda_1+ r(x)}{\alpha_2 + 1}  -  \frac{\lambda_2+ r(x)}{\alpha_1 + 1} }.
  \end{align}  
 \end{lemma}
 
 \begin{remark}\label{rem}
One can also consider other classes of conjugate priors  such as $\Gamma_1 = \{ \pi_{\alpha, \lambda_0}(\theta) :  \alpha\in [\alpha_1, \alpha_2], ~ \lambda_0   \text{ is fixed}\}$ or $
\Gamma_2  =  \{ \pi_{\alpha_0, \lambda}(\theta) :  \alpha=\alpha_0   \text{ is  fixed},  \lambda\in [\lambda_1, \lambda_2]  \}$.  The PRGM estimator of $\theta$ in $\Gamma_1$ or $\Gamma_2$ can be obtained using \eqref{gamma3} and  by letting $\lambda_1=\lambda_2=\lambda_0$ or $\alpha_1=\alpha_2=\alpha_0$, respectively. 
 \end{remark}

 \begin{example}\label{exp-prgm} In Example \ref{exponential-dist},  let  $\pi_{\alpha, \lambda}(\theta)\propto \theta^{\alpha-1} e^{-\theta\lambda}$ with   the posterior distribution $\pi(\theta| x) =\pi_{\alpha+1, \lambda+x}(\theta)$, and $\delta_{\pi}(x)= \frac{\alpha+1}{\lambda+x}$. Using \eqref{gamma3},  the  PRGM estimator of $\theta$ under the Stein loss function $L(\theta, \delta)= \frac{\delta}{\theta}- \log\frac{\delta}{\theta}-1$ in  $\Gamma = \{ \pi_{\alpha, \lambda}(\theta) ~: ~ \alpha\in [\alpha_1, \alpha_2] , ~ \lambda\in [\lambda_1, \lambda_2] \} $,   with $0<\alpha_1<\alpha_2$ and $0<\lambda_1<\lambda_2$ is given by
$$\delta^{\Gamma}_{PR}(X) = {\log \left ( \frac{\alpha_1+1}{\alpha_2+1} \frac{\lambda_1+X}{\lambda_2+X} \right)}\bigg/ \left({\frac{\lambda_1+X}{\alpha_2+1}-\frac{\lambda_2+X}{\alpha_1+1}}\right).$$
In $\Gamma_1$, as defined in Remark \ref{rem}, we have
$$\delta^{\Gamma_1}_{PR}(X)= \frac{(\alpha_1+1)\, (\alpha_2+1)}{\alpha_1 -\alpha_2}\log \left( \frac{\alpha_1+1}{\alpha_2+1}\right) \frac{1}{ \lambda_0 + X}.$$
Similarly, in  $\Gamma_2$,  we have 
$$\delta^{\Gamma_2}_{PR}(X)= \left(\frac{\alpha_0 + 1}{\lambda_2-\lambda_1}\right)\log \left(\frac{\lambda_2 + X}{\lambda_1 + X }\right).$$  \end{example}

%%%%%%%%%%%%%%%%%%%%%%%%%%%%
%%%%%%%%%%%%%%%%%%%%%%%%%%%
%%%%%%%%%%%%%%%%%%%%%%%%%%%

\section{\sc Intrinsic PRGM estimation }
In Section \ref{IPRGM-se}, we obtained the PRGM estimator of the natural parameter $\theta$  of the exponential family under the intrinsic loss function.  In some applications,  there may be interest in finding PRGM estimation of the original parameter of the underlying model rather than the natural parameter $\theta$. Unfortunately, like many other methods, PRGM estimators are not necessarily invariant under reparameterization.  
Although results of this nature,  that are not invariant under reparameterization,  can sometimes be interesting in theory,  they tend to be less useful in practice. Indeed, it is difficult to sell to a practitioner that the PRGM estimator of $h(\theta)$ is not necessarily  $h(\delta_{PR})$. In this section,  we obtain PRGM estimators that  are invariant under one-to-one smooth reparameterizations, hence the name intrinsic PRGM estimators.  

For  the exponential family \eqref{exp-family}, as opposed to the well known and commonly used conjugate prior \eqref{prior-density},   consider the following conjugate prior distribution  for $\theta$
\begin{align}\label{JCP}
\pi^J_{\alpha, \lambda}(\theta)  \propto  \{\beta(\theta)\}^{\alpha}\, e^{-\lambda \theta}\, \sqrt{I_{\theta}(\theta) },
\end{align}
where $I_{\theta}(\theta)$ is the Fisher information for $\theta$. 
Druilhet and Pommeret (2012)  introduced  \eqref{JCP} and referred to it as the \textit{Jeffrey's Conjugate Prior} (JCP). It is easy to see that the JCP is  invariant under  smooth reparameterizations, and   the necessary  conditions on $\alpha$ and $\lambda$ in \eqref{JCP},  leading to proper posterior distributions,  do not depend on the choice of the reparameterization.  The invariance property of JCP under  any smooth and one-to-one  reparameterization $\eta=h(\theta)$ can be shown  by the following relationship
$$I_{\eta}(\eta)= I_{\theta}(h^{-1}(\eta))\times  \big|\frac{d h^{-1}(\eta)}{d\eta} \big|^2. $$

\begin{remark} For the exponential family \eqref{exp-family}, since $I_{\theta}(\theta)= -H'(\theta)$, the JCP  is given by
$\pi^J_{\alpha, \lambda}(\theta)\propto \{\beta(\theta)\}^{\alpha} \, e^{-\lambda\theta}\sqrt{-H'(\theta)} $.
\end{remark}

First, we give the following result.

\begin{lemma}\label{invariant-bayes} Suppose $\delta^J_{\pi}$ is the Bayes estimator of the natural parameter $\theta$ of the exponential family \eqref{exp-family} under the intrinsic loss function \eqref{KL} with respect to the JCP distribution \eqref{JCP}. 
For every one-to-one smooth transformation $h(\theta)$,  the Bayes estimator of $h(\theta)$ is $h(\delta^J_{\pi})$. 
\end{lemma}

\noindent \textbf{Proof:}  The proof is similar to the proof of  Lemma 6.2 of Robert (1996) and hence  omitted. \\

%Let $\pi^J(\theta|x)$ be the posterior distribution of $\theta$ given $x$ with respect to the JCP distribution $\pi^J_{\alpha, \lambda}(\theta)$ as in \eqref{JCP}.  Suppose $\delta_{\pi}$ is the Bayes estimator of $\theta$ under the intrinsic loss function $L(\theta, \delta)= \log\left( \frac{\beta(\theta)}{\beta(\delta)}\right) + (\delta-\theta) \frac{\beta'(\theta)}{\beta(\theta)}$. That is 
%$$\delta_{\pi}(x)= argmin_{d}  \int _{\Theta}  L(\theta, d) \pi(\theta |x) d\theta.$$ 
% Using the invariance property of $\pi^J_{\alpha, \lambda}$, for any one-to-one transformation $\eta=h(\theta)$, $\pi^J(\theta| x) = \pi^J(\eta| x)$. Let $\tilde{\delta}_{\pi}$ be the Bayes estimator of $\eta$ under the loss $L(h(\theta), h(\delta)$. Then, 
% \begin{align*}
% \tilde{\delta}_{\pi} = 
% \end{align*}

Now, we state the main result of this section which can easily be  proved using the invariance property of both the class of JCP  distributions and the intrinsic loss functions under smooth reparameterization of $\theta$. 

\begin{theorem}\label{Invariant-PRGM} Suppose $\delta^{\Gamma^J}_{IPR}(X)$ is the PRGM estimator of the unknown parameter $\theta$ for the exponential family \eqref{exp-family} under the intrinsic loss function \eqref{KL-loss} with respect to a class $\Gamma^J$ of JCP distributions for $\theta$. Then,  for any one-to-one smooth transformation $h(\theta)$,   the PRGM estimator of $h(\theta)$ is $h(\delta^{\Gamma^J}_{IPR}(X))$. 
\end{theorem}
 
 \noindent\textbf{Proof:}
 By definition, the PRGM estimator of  $h(\theta)$ in the class $\Gamma^J$ of JCP distributions   is given by the solution of 
 $$\inf_{\delta}\sup_{\pi\in\Gamma^J} \rho(\delta^h_{\pi}, \delta)= \inf_{\delta}\sup_{\pi\in\Gamma^J} \left\{  \log\frac{\beta(\delta^h_{\pi})}{\beta(\delta)}+(\delta-\delta^h_{\pi})H(\delta^h_{\pi})\right\}, $$
where $\delta^h_{\pi}$ is the Bayes estimator of $h(\theta)$. 
 Note that   $\rho(\delta^h_{\pi}, \delta)= L(\delta^h_{\pi}, \delta)$ where $L$ is defined in \eqref{KL-loss}. Now, using the invariance property  of $L$ and Lemma \ref{invariant-bayes},  since $\delta^h_{\pi}= h(\delta_{\pi})$, with  $\delta_{\pi}$ being  the Bayes estimator of $\theta$, we have 
 \begin{align*}
 \inf_{\delta} \sup_{\pi\in\Gamma^J}\rho  (\delta^h_{\pi}, \delta)
 &=  \inf_{\delta} \sup_{\pi\in\Gamma^J}\rho  (h(\delta_{\pi}), \delta)\\
 &=  \inf_{t: h(t)= \delta } \sup_{\pi\in\Gamma^J}\rho  ( h(\delta_{\pi}), h(t))\\
 &= \inf_{t: h(t)= \delta } \sup_{\pi\in\Gamma^J}\rho  (\delta_{\pi}, t).
 \end{align*}
 Therefore, if $\delta^{\Gamma^J}_{IPR}(X)$  is the PRGM estimator of $\theta$, i.e., $\delta^{\Gamma^J}_{IPR}$  minimizes  (in $t$)  $\sup_{\pi\in\Gamma^J}\rho  (\delta_{\pi}, t)$, then,  the transform $h(\delta^{\Gamma^J}_{IPR}(X))$  is the PRGM estimator of $h(\theta)$, that is, $h(\delta^{\Gamma^J}_{IPR})$ minimizes  (in $\delta$)  $\sup_{\pi\in\Gamma^J}\rho  (\delta^h_{\pi}, \delta)$ and this completes the proof. $\Box$

 \begin{example}\label{ex-exp-bayes} Suppose $X\sim Exp(\sigma)$ with $\sigma, x>0$. In Example \ref{exponential-dist},  we showed that  the intrinsic loss  for estimating $\theta=\sigma^{-1}$ by $\delta$  reduces to   the Stein  loss function
$$L(\theta, \delta) =\frac{\delta}{\theta} - \log \frac{\delta}{\theta} -1. $$
Under the  JCP distribution  $\pi^J_{\alpha, \lambda}(\theta) \propto \theta^{\alpha-2} e^{-\theta\lambda}$,  $\alpha>1$,  the posterior distribution is a $Gamma(\alpha, \frac{1}{\lambda+x})$ with $\pi^J(\theta| x)\propto \theta^{\alpha-1} e^{-(\lambda+ x) \theta}$ which results in the  Bayes estimator of $\theta$  as  $\delta_{\pi}(X)= \frac{\alpha}{\lambda+X}$. Also,  the intrinsic  PRGM estimator of $\theta$ under $L(\theta, \delta)$ is given by
$$\delta^{\Gamma^J}_{IPR}(X)=\frac{\log\frac{1}{\overline{\delta}(X)} - \log\frac{1}{\underline{\delta}(X)}}{\frac{1}{\overline{\delta}(X)}-\frac{1}{~\underline{\delta}(X)}}.$$
Now, for the  estimation of  $\eta=\sigma=\frac{1}{\theta}$ using $\tilde{\delta}$, it is easy to see that   the Bayes estimator of $\eta$ under the entropy  loss function $$L(\eta, \tilde{\delta})= \frac{\eta}{\tilde{\delta}}- \log  \frac{\eta}{\tilde{\delta}}-1, $$
is given by $ \tilde{\delta}_{\pi}(X)= \frac{\lambda+X}{\alpha}= \frac{1}{\delta_{\pi}(X)}$. To see this, note that $\pi^J(\eta)\propto \eta^{-\alpha} e^{-\lambda/ \eta}$ with $\pi^J(\eta| x)\propto \eta ^{-(\alpha+1)} e^{-\frac{\lambda+x}{\eta}}$ and $\tilde{\delta}_{\pi}(x)= E[\eta|x]$. Also, the intrinsic  PRGM estimator of $\eta$  is given by 
$$\tilde{\delta}^{\Gamma^J}_{IPR}(X) = \frac{\overline{\tilde{\delta}}(X) - \underline{\tilde{\delta}}(X)}{\log\overline{\tilde{\delta}}(X) - \log\underline{\tilde{\delta}}(X)}= \frac{\frac{1}{\overline{\delta}(X)}-\frac{1}{~\underline{\delta}(X)}}{\log\frac{1}{\overline{\delta}(X)} - \log\frac{1}{\underline{\delta}(X)}}= \frac{1}{\delta^{\Gamma^J}_{IPR}(X)}. $$
For the   PRGM estimation  of $\theta$ under the Entropy loss function and its application to record data analysis we refer to  Jafari Jozani and Parsian (2008). Similarly, if $\eta^*=-\frac{1}{a} \log \theta$, $\alpha\neq 0$, then the intrinsic  PRGM estimator of $\eta^*$ under the LINEX loss function 
$$L(\eta^*, \delta^*) = e^{a\, (\eta^*-\delta^*) }- a\, (\eta^*-\delta^*)-1,  $$
is given by 
$${\delta^*}^{\Gamma^J}_{IPR}(X)= \underline{\delta^*}(X) + \frac 1a \log\left\{ \frac{e^{a\, (\overline{\delta^*}(X) - \underline{\delta^*}(X)) }-1}{a\, (\overline{\delta^*}(X) - \underline{\delta^*}(X))}\right\}= \frac{1}{a} \log \delta ^{\Gamma^J}_{IPR}(X), $$
which is the PRGM estimator obtained in Boraty\'nska (2006).  
\end{example}

For the exponential family \eqref{exp-family},  suppose that the prior distribution  belongs to the following class of JCP distributions: 
\begin{align}
\Gamma^J & = \{ \pi^J_{\alpha, \lambda}(\theta) ~: ~ \alpha\in [\alpha_1, \alpha_2], ~ \lambda\in [\lambda_1, \lambda_2]\}, 
\end{align}  
 for suitable choices  of $\alpha_1<\alpha_2$ and $\lambda_1<\lambda_2$.   
We continue  with some  applications of  Theorem \ref{Invariant-PRGM} under the above class of priors.  Similar results can  be obtained in other classes of JCP distributions (see Remark \ref{rem}),  which we do not present here. In view of Theorem \ref{Invariant-PRGM}, and to obtain an intrinsic  PRGM estimator,  the critical condition  is that the elements of the underlying  class of  prior distributions  are in the form of \eqref{JCP} and the underlying loss function   is intrinsic. We observe that, in many cases (see  Examples \ref{exp-prgm-invariant} and \ref{binom-prgm-invariant})  intrinsic   PRGM estimators under $\Gamma^J$ can be obtained using the PRGM estimators  under  the usual class $\Gamma$  of conjugate priors with  modified   values of $\alpha_i$s and $\lambda_i$s in $\Gamma$, $i=1, 2$.  One can easily check that this will happen whenever the mean-value parameter is conjugate for the natural parameter in the sense of Gutierrez-Pena and Smith (1995). In the one-parameter case, a sufficient condition for this is that the exponential family have a quadratic variance function (see Section 3.3 of Gutierrez-Pena and Smith (1995)).

 \begin{example}\label{exp-prgm-invariant} In Example \ref{ex-exp-bayes},  we showed that $\pi_{\alpha, \lambda}^J(\theta)\propto \theta^{\alpha-2} e^{-\theta\lambda}$ and $\delta_{\pi}(x)= \frac{\alpha}{\lambda+x}$. Since $\pi_{\alpha, \beta}^J(\theta|x)$ is equal to $\pi(\theta| x)$,  the posterior distribution of $\theta$,  given the usual conjugate prior  $ \pi_{\alpha-1, \lambda+x}(\theta)$, the intrinsic PRGM  estimator of $\theta$ under the Stein loss function and the class of JCP distributions   can be obtained using the PRGM estimator of $\theta$  under the usual class of conjugate priors ( as in Example \ref{exp-prgm}),  by replacing  $\alpha_i$ with $\alpha_i-1$, $i=1, 2$. For example,   the intrinsic  PRGM estimator of $\theta$ in  $\Gamma^J$  with $0<\alpha_1<\alpha_2$ and $0<\lambda_1<\lambda_2$ is given by 
 $$\delta^{\Gamma^J}_{IPR}(X) = {\log \left ( \frac{\alpha_1}{\alpha_2} \frac{\lambda_1+X}{\lambda_2+X} \right)}\bigg/ \left({\frac{\lambda_1+X}{\alpha_2}-\frac{\lambda_2+X}{\alpha_1}}\right).$$
 Let $\Gamma_1^J=\{ \pi_{\alpha, \beta}^J(\theta): \alpha\in[\alpha_1, \alpha_2]  \text{ and } \lambda= \lambda_0\} $. Then,  
   the intrinsic  PRGM estimator of $\theta$ in  $\Gamma _1^J$,   with $0<\alpha_1<\alpha_2$ and $\lambda_0>0$,  is given by
     $$\delta^{\Gamma_1^J}_{IPR}(X)=\left( \frac{\alpha_1\, \alpha_2}{\alpha_1 -\alpha_2}\right)\log \left( \frac{\alpha_1}{\alpha_2}\right) \frac{1}{ \lambda_0 + X}.$$
Similarly, in $\Gamma_2^J=\{ \pi_{\alpha, \beta}^J(\theta): \alpha=\alpha_0 \text{ fixed } \text{ and } \lambda\in[\lambda_1, \lambda_2]\} $, 
   the intrinsic PRGM estimator of $\theta$ in  $\Gamma _2^J$,   with $\alpha_0>0$ and $\lambda_1, \lambda_2>0$,  is given by
   $$\delta^{\Gamma^J_2}_{IPR}(X)= \frac{\alpha_0}{\lambda_2-\lambda_1}\log \left\{ \frac{\lambda_2 + X}{\lambda_1 + X}\right\}.$$ Similar results can be obtained for estimating any smooth and one-to-one   function   of $\theta$  under corresponding class of JCP  distributions.  
%Now,  using Theorem \ref{Invariant-PRGM} the PRGM estimator of $\eta=\frac{1}{\theta}$ under $L(\eta, \tilde{\delta})= \frac{\eta}{\tilde{\delta}}-\log \frac{\eta}{\tilde{\delta}} -1$ under $\Gamma_i^J$, $i=1, 2, 3$,  is simply given by $\tilde{\delta}_{IPR, i}(x)= \frac{1}{\delta_{IPR, i}(x)}$.
  \end{example}

 \begin{example}\label{binom-prgm-invariant} (Binomial Distribution). In Example \ref{binomial-dist}, we showed that pmf of $X$ can be written as $f(x| \theta)= \binom{n}{x} (\frac{e^{\theta}}{1+ e^{\theta}})^n e^{-x\theta} $ with $\theta=\log(\frac{1-p}{p})$. Here $I_{\theta}(\theta) = \frac{e^{\theta}}{(1+ e^{\theta})^2}$ and the JCP for $\theta$ is obtained as $\pi^J_{\alpha, \lambda}(\theta)\propto(\frac{e^{\theta}}{1+ e^{\theta}})^{\alpha} e^{- \lambda\theta}\frac{e^{\theta/2}}{1+ e^{\theta}}$. This results in the posterior distribution 
  $\pi^J(\theta|x)\propto(\frac{e^{\theta}}{1+ e^{\theta}})^{\alpha+n+1} e^{-(x+ \lambda+ \frac12)\theta}. $ Since $\pi^J(\theta|x)$ is equal to $\pi(\theta| x)$, the  posterior distribution of $\theta$, given   the usual conjugate prior  $\pi_{\alpha+1, \lambda+ \frac12}(\theta)$,  the intrinsic  PRGM estimator $\delta^{\Gamma^J}_{IPR}(X)$ of  $\theta$ in  $\Gamma^J$ can be obtained  using   \eqref{prgm-binom} and  by replacing $\alpha_i$ and $\lambda_i$ with $\alpha_i+1 $ and $\lambda_i+ \frac12 $, $i=1, 2$,  respectively.  
 Also, the intrinsic  PRGM estimator of $p= \frac{1}{1+ e^{\theta}}$ under the  loss function 
 \begin{align*}
 L(p, \tilde{\delta}) = p\, \log \left(\frac{p}{\tilde{\delta}}\right) + (1-p)\, \log \left( \frac{1-p}{1-\tilde{\delta}}\right),
 \end{align*} 
 is given by $\tilde{\delta}^{\Gamma^J}_{IPR}(X)=\{1+ e^{\delta_{IPR}^{\Gamma^J}(X)}\}^{-1}$.
 \end{example}
 
%  \begin{remark} Note that similar development,   as in Section 2.1,   can be obtained concerning the Bayesianity of %the invariant PRGM estimators with respect to a prior in the class of JCP distributions. 
  %\end{remark}
%
%
%\begin{example}
%\label{boundednormal}
%For a normal model $X \sim \hbox{N}(\theta; \sigma^2)$ with $|\theta| \leq m$ and known $\sigma^2$, Rios Insua, Ruggeri and Vidakovic (1995) obtain
%the following  posterior regret $\Gamma$-minimax solutions  under squared-error loss:
%\begin{enumerate}
%\item[ (a)] $\delta^{PRGM}(x) = \frac{m}{2} \tanh(\frac{mx}{\sigma^2})$, for the class $\Gamma$ of symmetric priors on $[-m,m]$;
%
%\item[ (b)] $\delta^{PRGM}(x) = \frac{1}{2} \delta_U(x)$, for the class $\Gamma$ of unimodal, symmetric, absolutely continuous
%priors on $[-m,m]$, where $\delta_U$ is the Bayes estimator with respect to a uniform prior on $[-m,m]$ given by
%$\delta_u(x) = x + \sigma \left([\phi(\frac{x+m}{\sigma}) - \phi(\frac{x-m}{\sigma})]/[\Phi(\frac{x+m}{\sigma}) - \Phi(\frac{x-m}{\sigma})]\right)$;
%with $\phi$ and $\Phi$ representing the pdf and cdf of the standard normal.
%\end{enumerate}
%
%\end{example}
%

%%%%%%%%%%%%%%%%%%%%%%%%%%%%%%%
%%%%%%%%%%%%%%%%%%%%%%%%%%%%%%%
%%%%%%%%%%%%%%%%%%%%%%%%%%%%%%%

\section {\sc  PRGM, Intrinsic PRGM  and Bayes estimators  }

In this section,  we provide some general results concerning the Bayesianity of  the PRGM and intrinsic PRGM estimators of $\theta$ for the exponential family distribution \eqref{exp-family}   under the  intrinsic loss function \eqref{KL-loss}  with respect to  priors in the underlying class of prior distributions. The results are  only  presented  for PRGM estimators of $\theta$, but they can also be used  for  intrinsic PRGM estimators  by  simple modifications.   Our framework in this section closely resembles the one introduced by Rios Insua et al.\ (1995), who considered  similar problem for  the quadratic loss function.  Results of  this nature are also obtained by Zen and DasGupta (1993) under the quadratic loss function for the  binomial distribution.     Several of the following  preliminary results and detailed proofs are reported here for the sake of completeness.  The idea is to check  the continuity of  the underlying Bayes estimator with respect to the prior. Similar to Rios Insua et al.\ (1995) we study two cases, when (a) the class of prior distributions is convex, or (b) the underlying class of prior distributions depends on a hyper-parameter belonging to a connected set.   

 First,  consider the  situation where the  class $\Gamma$  of priors is  convex. That is,  if $\pi_{0}, \pi_{1} \in\Gamma$,  then  
$\pi_{t}=t\pi_{0}+(1-t)\pi_{1}$ belongs to $\Gamma$, for any $t\in[0, 1]$. Suppose that  $X$ is a
random variable whose  density belongs to the family of distributions \eqref{exp-family}. 
%We showed that the  Bayes estimator of $\theta$ under the intrinsic loss function \eqref{KL-loss} is given as the solution of 
%$H(\delta_{\pi}(x))=E(\frac{\beta'(\theta)}{\beta(\theta)}| x)$.
Let
$\psi(t)=H(\delta_{\pi_{t}}(x))$
which is a decreasing function of $\delta_{\pi_{t}}$ for any $t\in[0,1]$. In
the next lemma we show that $\psi(t)$ is a  continuous function in its domain   $t\in [0,1]$.

\begin{lemma}\label{continuous}
Suppose  $\psi(t)$,   the posterior expectation of
$H(\theta)=\frac{\beta'(\theta) }{\beta(\theta)}$  when 
$\pi_{t}=t\pi_{0}+(1-t)\pi_{1}$, $t\in[0,1]$, is finite. Then,  $\psi(t)$ is continuous in $t$, 
 $t\in[0,1]$.
\end{lemma}
\noindent\textbf{Proof.} Let $a_{i}=
\int_{\Theta}H(\theta)\pi_{i}(\theta)f(x|\theta)d\theta$,  $m_{i}(x)=\int_{\Theta}\pi_{i}(\theta)f(x|\theta)d\theta$ for
$i=1,2$, and suppose that $a_i$s and $m_i$s exist and are finite. Then
\begin{align*}
\psi(t) 
&= E_{\pi_t}[H(\theta)|x] \cr
%&=\int_{\Theta}H(\theta)\times\frac{f(x|\theta)\pi_{t}(\theta)}{m_{t}(x)}d\theta\cr
%&=\int_{\Theta}\frac{\beta'(\theta)}{\beta(\theta)}.\frac{f(x|\theta)(t\pi_{0}(\theta)+(1-t)\pi_{1}(\theta))}
%{\int_{\Theta}\pi_{t}(\theta)f(x|\theta)d\theta}d\theta\cr
&=\frac{ t\int_{\Theta}H(\theta)\, f(x|\theta)\, \pi_{0}(\theta)d\theta+(1-t)\int_{\Theta}H(\theta)\, 
f(x|\theta)\, \pi_{1}(\theta)d\theta}{ t\, \int_{\Theta}f(x|\theta)\pi_{0}(\theta)d\theta + (1-t)\int_{\Theta}f(x|\theta)\pi_{1}(\theta)d\theta}\cr
&=\frac{t a_{0}+(1-t) a_{1}}{t m_{0}(x)+(1-t)m_{1}(x)},
\end{align*}
which is a continuous function of $t$, $t\in[0,1]$.  $~\Box$
 
 Now, we use the continuity of $\psi(t)$ to prove that,  under the conditions of Lemma \ref{continuous},   the PRGM estimator 
$\delta_{PR}$ is  Bayes if the class of priors is convex.

\begin{theorem}\label{convex-class}  Suppose $\Gamma$ is a convex class of prior distributions on the unknown parameter $\theta$ of the exponential family of distributions \eqref{exp-family}. Then, there exists a prior distribution $\pi\in\Gamma$ such that $\delta_{PR}= \delta_{\pi}$, where $\delta_{PR}$ is defined in \eqref{PRGM-est}. 
\end{theorem}
\noindent\textbf{Proof.}  Following the definition  of $\underline{\delta}$ and $\overline{\delta}$, consider a small enough  $\varepsilon>0$ and two    prior distributions 
$\pi_{0}, \pi_{1}\in\Gamma$ such that
$$\delta_{\pi_{0}}<\underline{\delta}+\epsilon<\delta_{PR}<\overline{\delta}-\epsilon<\delta_{\pi_{1}}.$$
Since $H(\cdot) $ is a decreasing 
function, then
\begin{align*}
H(\delta_{\pi_1})<H(\overline{\delta}-\epsilon)<H(\delta_{PR})<H(\underline{\delta}+\epsilon)<
H(\delta_{\pi_{0}}).
\end{align*}
Let $\pi_{t}=t\pi_{0}+(1-t)\pi_{1}$, $t\in[0,1]$  and define 
$\psi(t)=H(\delta_{\pi_{t}})$. Note that,  
$\psi(0)=H(\delta_{\pi_{0}})$ 
and
$\psi(1)=H(\delta_{\pi_{1}})$.
Now, from Lemma \ref{continuous}, the continuity of $\psi(t)$ in $t$ shows  that there exists a $t^*\in[0,1]$ such that
$
\psi(t^*)=H(\delta_{\pi_{t^*}})=H(\delta_{PR}),
$
which completes the proof. $\Box$\\

A shortcoming of the  result in Theorem \ref{convex-class} is that it is not applicable to the cases where the class of prior distributions depends on a hyper-parameter whose range is connected. For this case, we prove Lemma \ref{rios1}  and Theorem \ref{rios2}  which are   simple extensions of  Lemma 3.2 and Proposition 3.2 of Rios Insua et al.\ (1995).  The proof of Lemma \ref{rios1} is essentially similar to the proof of Lemma 3.2 of Rios Insua et al.\ (1995). The same is true of Theorem \ref{rios2}. Nonetheless,  we provide the proofs in the Appendix for  the sake of completeness. 
Let 
\begin{eqnarray}
\psi(\pi)
=\frac{\int_{\Theta}H(\theta)f(x|\theta)\pi(\theta )d\theta}{\int_{\Theta}f(x|\theta)\pi(\theta)d\theta}
=\frac{r(\pi)}{s(\pi)}. 
\end{eqnarray}
Consider $d(\pi,\pi')=\sup_{\Theta}|\pi(\theta)-\pi'(\theta)|$  to be the usual  $l_{\infty }$  distance between prior  densities $\pi$ and $\pi'$,  where  $H(t)= \beta'(t)/\beta(t)$  is defined as before. 

\begin{lemma}\label{rios1}
Suppose that  $\int_{\Theta}\big |H(\theta)\big
|f(x|\theta)d\theta$ exists and it is finite. Then, $\psi(\pi)$ is continuous 
in $\pi$, in the topology generated by the $l_{\infty}$ distance.
\end{lemma}

\begin{theorem}\label{rios2} Let
$\Gamma=\{\pi_{\alpha}:\alpha\in\Lambda\}$,  where $\Lambda$ is a
connected set and $\pi_{\alpha}$'s are densities. Under the conditions of Lemma \ref{rios1} and the assumption that  
$\alpha_{n}\rightarrow\alpha$ implies 
$d(\pi_{\alpha_{n}},\pi_{\alpha})\rightarrow 0$,  there exists a prior distribution 
$\pi\in\Gamma$ such that $\delta_{PR}=\delta_{\pi}$, that is, the PRGM estimator \eqref{PRGM-est} is Bayes.
\end{theorem}

 %In the following example,  we consider a case where the PRGM  (or intrinsic PRGM) estimator  is Bayes with %respect to the same prior in the underlying class of prior distribution, regardless of the observed value of $x$. 

  In the following lemma,  we provide a sufficient condition under which  the PRGM  (or intrinsic PRGM) estimator  is Bayes with respect to the same prior in the underlying class of prior distribution, regardless of the observed value of $x$. 
 
\begin{lemma}
Let $\Gamma=\{ \pi_{\alpha}: \alpha\in[\alpha_1, \alpha_2]\}$ be the class of prior distributions. Suppose    the Bayes estimator $\Psi(\alpha, x)= H^{-1}\{E[H(\theta) | x) \}$  is a differentiable  function of the hyper-parameter $\alpha$ and the observed value $x$. Assume that we are  under the conditions of Theorem \ref{rios2}.  If 
\begin{align}\label{prior-ind}\frac{\partial}{\partial x}\Psi(\alpha, x)= \frac{\partial}{\partial x} \left\{\frac{\Psi(\alpha_1, x)\, H(\Psi(\alpha_1, x))-\Psi(\alpha_2, x)\, H(\Psi(\alpha_2, x))-\log\frac{\beta(\Psi(\alpha_1, x))}{\beta(\Psi(\alpha_2, x))}}{ H(\Psi(\alpha_1, x))-H(\Psi(\alpha_2, x))}
  \right\},  \end{align}
has a constant solution in $\alpha$, then there is a data independent prior $\pi_{\alpha}\in\Gamma$ resulting  in the PRGM estimate as the Bayes estimate of  the natural parameter $\theta$ of the exponential family \eqref{exp-family} under the intrinsic loss function \eqref{KL}. 
\end{lemma}
 \noindent\textbf{Proof:} Under the conditions of Theorem \ref{rios2}, there exists a solution  $\alpha(x)$ such that the PRGM estimator \eqref{PRGM-est} is  Bayes with respect to the prior $\pi_{\alpha(x)}\in \Gamma$ under the intrinsic loss function \eqref{KL}. That is, 
 $$\Psi(\alpha(x), x)= \frac{\Psi(\alpha_1, x)\, H(\Psi(\alpha_1, x))-\Psi(\alpha_2, x)\, H(\Psi(\alpha_2, x))-\log\frac{\beta(\Psi(\alpha_1, x))}{\beta(\Psi(\alpha_2, x))}}{ H(\Psi(\alpha_1, x))-H(\Psi(\alpha_2, x))}.   $$
Now, differentiating the equation with respect to $x$ leads to 
\begin{align*} &\frac{\partial}{\partial \alpha} \Psi(\alpha(x), x) \frac{d\alpha(x)}{dx} + \frac{\partial}{\partial x} \Psi(\alpha(x), x)\\
&= \frac{\partial}{\partial x} \left\{\frac{\Psi(\alpha_1, x)\, H(\Psi(\alpha_1, x))-\Psi(\alpha_2, x)\, H(\Psi(\alpha_2, x))-\log\frac{\beta(\Psi(\alpha_1, x))}{\beta(\Psi(\alpha_2, x))}}{ H(\Psi(\alpha_1, x))-H(\Psi(\alpha_2, x))}
  \right\}.   \end{align*}
If $\alpha(x)$ is data independent, i.e., $\alpha(x)= \alpha$, then $\frac{d\alpha(x)}{dx}=0$. Now, the desired value for    $\alpha$ is the constant solution to the equation \eqref{prior-ind} leading to a data independent prior for the PRGM estimator to be Bayes.  $\Box$

\begin{example}
In Example \ref{normal}, the condition \eqref{prior-ind} reduces to the condition (5) in Proposition 3.3 of Rios Insua et al.\ (1995) as follows
$$2 \frac{\partial}{\partial x}\Psi(\alpha, x) = \frac{\partial}{\partial x}\Psi(\alpha_1, x)+ \frac{\partial}{\partial x}\Psi(\alpha_2, x).$$
Now, consider the class $\Gamma=\{ \pi_{\alpha, \lambda_0}: \alpha\in[\alpha_1, \alpha_2], \lambda_0 \text{ is fixed }\}$ of conjugate priors where $\pi_{\alpha, \lambda_0}$ is given by \eqref{prior-density} with $\theta=\mu$ and $\beta(\theta) = e^{-\theta^2/2}$. Here,  the Bayes estimator of $\theta$ is given by $\Psi(\alpha, X)=\delta_{\pi_{\alpha, \lambda}}(X)= \frac{X-\lambda_0}{\alpha+1}$. It is easy to see that,  the PRGM estimator of $\theta$  given by 
$$\delta_{PR}(X)= \frac12 \left\{ \frac{X-\lambda_0}{\alpha_1+ 1} + \frac{X-\lambda_0}{\alpha_2+1} \right\}, $$
is Bayes with respect to the data independent prior $\pi_{\alpha^*, \lambda_0} \in \Gamma$ where $\alpha^*$ is given as the solution to the following equation
$$ \frac{2}{ \alpha^*+1}= \frac{1}{\alpha_1+1} + \frac{1}{\alpha_2+1}. $$
That is, $\alpha^*= \frac{\alpha_1+\alpha_2+ 2\alpha_1\alpha_2}{\alpha_1+ \alpha_2 + 2}\in [\alpha_1, \alpha_2]$ and $\delta_{PR}(X)= \frac{X-\lambda_0}{\alpha^* +1}= \delta_{\pi_{\alpha^*, \lambda}}(X)$. 

\end{example}

\begin{example}
In Example \ref{exp-prgm}, the condition \eqref{prior-ind} reduces to
$$\frac{\partial}{\partial x}\Psi(\alpha, x)= \frac{\partial}{\partial x} \left\{ \frac {\log \frac{1}{\Psi(\alpha_1, x)} - \log \frac{1}{\Psi(\alpha_2, x)}}  {\frac{1}{\Psi(\alpha_1, x)} - \frac{1}{\Psi(\alpha_2, x)}} \right\}.  $$
Now, consider the class $\Gamma_1 = \{ \pi_{\alpha, \lambda_0}(\theta) :  \alpha\in [\alpha_1, \alpha_2], ~ \lambda_0   \text{ is fixed}\}$  of conjugate priors on $\theta$. Here, the Bayes estimator of $\theta$ with respect to the prior $\pi_{\alpha, \lambda_0}(\theta)$  is $\Psi(\alpha, X)= \delta_{\pi_{\alpha, \lambda_0}}(X)= \frac{\alpha+1}{\lambda_0 + X}$. The PRGM estimator of $\theta$ is then  Bayes with respect to a data independent prior $\pi_{\alpha^*, \lambda_0}(\theta)\in\Gamma_1$, if  there exists  a data independent  solution $\alpha^*$ to the equation  
$$-\frac{\alpha^*+ 1}{(\lambda_0+X)^2} =  - \log \left( \frac{\alpha_1+1}{\alpha_2+1}\right)\frac{(\alpha_1+1)\, (\alpha_2+1)}{\alpha_1 -\alpha_2} \frac{1}{(\lambda_0+X)^2}.   $$
A straightforward calculation shows that 
$$ \alpha^*= \frac{(\alpha_1+1)\, (\alpha_2+1)}{\alpha_1 -\alpha_2}\log \left( \frac{\alpha_1+1}{\alpha_2+1}\right)-1 \in [\alpha_1, \alpha_2].$$ Therefore, 
the PRGM estimator of $\theta$ under the Stein loss function can be obtained as the Bayes estimator of $\theta$ with respect to the  prior distribution $\pi_{\alpha^*, \lambda_0}(\theta)\in\Gamma_1$ as follows
$$\delta^{\Gamma_1}_{PR}(X)= \frac{(\alpha_1+1)\, (\alpha_2+1)}{\alpha_1 -\alpha_2}\log \left( \frac{\alpha_1+1}{\alpha_2+1}\right) \frac{1}{ \lambda_0 + X} =\frac{\alpha^*+1}{\lambda_0 + X}=\delta_{\pi_{\alpha^*, \lambda_0}}(X).  $$ 
 Similarly, in Example \ref{exp-prgm-invariant}, one can easily show that the intrinsic PRGM estimator $\delta^{\Gamma_1^J}_{IPR}(X)$ is the Bayes estimator of $\theta$ under the Stein loss function with respect to the prior distribution $\pi^J_{\alpha^{**}, \lambda_0}\in \Gamma^J_1$, when $\alpha^{**}=  \frac{\alpha_1\, \alpha_2}{\alpha_1 -\alpha_2}\log ( \frac{\alpha_1}{\alpha_2})$. Note that ${1}/{\alpha^{**}}$ is the logarithmic mean of $1 / \alpha
_1$ and $1/\alpha_2$, and  $\alpha^{**}\in[\alpha_1, \alpha_2]$.\end{example}

%%%%%%%%%%%%%%%%%%%%%%%%%%%%%%%%%
%%%%%%%%%%%%%%%%%%%%%%%%%%%%%%%%%
%%%%%%%%%%%%%%%%%%%%%%%%%%%%%%%%%

\section{\sc Concluding Remarks}
Invariant estimators  are  usually  demanding  in practice.  In this paper,  we have provided general results  concerning the PRGM estimation of the natural parameter of the one-parameter exponential family of distributions under intrinsic loss functions. The PRGM estimators are shown to be invariant  to one-to-one smooth reparameterizations under intrinsic loss functions and the class  of Jeffrey's conjugate prior distributions. Moreover, when the class of priors are convex or  dependant on a hyper-parameter belonging to  a connected set, we show that the obtained  PRGM estimators could be Bayes with respect to  prior distributions in the underlying class of priors. Several examples  are provided to clarify the results. 

 \section*{\sc Acknowledgements}
Mohammad Jafari Jozani gratefully  acknowledges  the partial support of  the Natural Sciences and  Engineering Research Council of Canada.  This work was done during the second author's visit  to the University of Manitoba, Department of Statistics. 

\noindent 

\section{\sc Appendix}
\subsection{Proof of Lemma \ref{rios1}}
  Suppose that
$d(\pi,\pi')<\epsilon$. Then, for all $\theta\in \Theta$, 
$\pi(\theta)-\epsilon\leq\pi'(\theta)\leq\pi(\theta)+\epsilon, 
$
and so 
\begin{align}\label{fpi}
f(x|\theta)\pi(\theta)-f(x|\theta)\epsilon\leq
f(x|\theta)\pi'(\theta)\leq
f(x|\theta)\pi(\theta)+f(x|\theta)\epsilon.
\end{align}
Upon  integrating \eqref{fpi} over $\theta$ we get
\begin{align*}
s(\pi)-\epsilon\int_{\Theta}f(x|\theta)d\theta\leq s(\pi')\leq
s(\pi)+\epsilon\int_{\Theta}f(x|\theta)d\theta.
\end{align*}
Let $\theta_0\in\Theta$ be such that $H(\theta)=\frac{\beta'(\theta)}{\beta(\theta)} >0$ for all $\theta <\theta_0$,  and  $H(\theta)\leq 0$ for all $\theta\geq \theta_0$. For $\theta< \theta_0$, 
multiplying \eqref{fpi} by
$H(\theta)\geq0$ results in 
\begin{align}\label{+bpi}
H(\theta)f(x|\theta)\pi(\theta)-\epsilon\, H(\theta)f(x|\theta)\leq
H(\theta)f(x|\theta)\pi'(\theta)\leq
H(\theta)f(x|\theta)\pi(\theta)+\epsilon\, H(\theta)f(x|\theta),
\end{align}
while for $\theta\geq\theta_0$ we  have
\begin{align}\label{-bpi}
H(\theta)f(x|\theta)\pi(\theta)+\epsilon\, H(\theta)f(x|\theta)\leq
H(\theta)f(x|\theta)\pi'(\theta)\leq
H(\theta)f(x|\theta)\pi(\theta)-\epsilon\, H(\theta)f(x|\theta).
\end{align}
Using   \eqref{+bpi} and \eqref{-bpi} and integrating over $\theta$, leads to 
\begin{align*}\label{abpi}
r(\theta)-\epsilon\int_{\Theta}\big
|H(\theta)\big |f(x|\theta)d\theta\leq
r'(\theta)\leq r(\theta)+\epsilon\int_{\Theta}\big
|H(\theta)\big |f(x|\theta)d\theta.
\end{align*}  Since we assumed  that $\int_{\Theta}\big
|H(\theta)\big
|f(x|\theta)d\theta=K_{1}<\infty$, then $\int_{\Theta}
f(x|\theta)d\theta=K_{2}<\infty$,  and 
\begin{align*}
\frac{r(\pi)-\epsilon K_{1}}{s(\pi)+\epsilon
K_{2}}\leq\frac{r(\pi')}{s(\pi')}\leq\frac{r(\pi)-\epsilon
K_{1}}{s(\pi)-\epsilon K_{2}}.
\end{align*}
If $\epsilon\rightarrow 0$,  then $\frac{r(\pi')}{s(\pi')}\rightarrow
\frac{r(\pi)}{s(\pi)}$. Therefore, $\psi(\pi')\rightarrow\psi(\pi)$.

%%%%%%%%%%%%%%%%%%%%%%%%%%
\subsection{Proof of Theorem \ref{rios2}}
 We consider
$\pi_{0}=\pi_{\alpha_{0}}$ and $\pi_{1}=\pi_{\alpha_{1}}$. Due to the 
connectedness of $\Lambda$,  there is a continuous path  $g(t)\in\Lambda, t\in[0,1]$ such
that $g(0)=\alpha_{0}$,  and $g(1)=\alpha_{1}$. Let
$\psi(t)=\psi(\pi_{g(t)})=H(\delta_{\pi_{g(t)}}),
t\in[0,1]$.  $\psi(t)$ is  a continuous function  in $t$, so there is
$t^*\in[0,1]$ such that $\psi(t^*)=\psi(\pi_{PR})$ leading to $\pi_{g(t^*)}\in\Gamma$ as   the  prior distribution  we were
looking for.

% Hence
%$$\frac{\beta'(\delta_{\pi_{g(t')}})}{\beta(\delta_{\pi_{g(t')}})}=\frac{\beta'(\delta_{\pi}^{P})}{\beta(\delta_{\pi}^{P})}.$$
%Therefore $\delta_{\pi_{g(t')}}=\delta_{\pi}^{P}$
%

\subsection*{\sc References}
\small
\renewcommand{\baselinestretch}{1.0}
\begin{description}

%
%\item \vskip  -2mm Berger, J.O. (1985). {\it Statistical Decision Theory and Bayesian Analysis.} New York, Springer-Verlag.
%
\item \vskip  -2mm Berger, J.O. (1994). {\it An overview of robust Bayesian analysis}. Test, \textbf{3}, 5--124.

\item\vskip-2mm Bernardo, J.M. (2011). {\it Integrated objective Bayesian estimation and hypothesis testing (with discussion)}. In Bayesian Analysis 9 (eds. J.M. Bernardo, M.J. Bayarri, J.O. Berger, A.P. Dawid,  D. Heckerman, A.F.M. Smith and M. West). Oxford University Press, 1--68. 

\item\vskip-2mm Bernardo, J.M. and Smith, A.F.M. (1994). {\it Bayesian Theory}. Chichester: Wiley.

%\item\vskip -2mm Betro, B., and Ruggeri, F. (1992). {\it
%Conditional $\Gamma$-minimax actions under convex losses}.
%Communications in Statistics: Theory and Methods, {\bf 21},
%1051-1066.

%\item\vskip -2mm Boraty\'{n}ska, A. (1997). {\it Stability of Bayesian inference in exponential families}.
%Statistics \& Probability Letters, {\bf 36}, 173--178.

%\item \vskip  -2mm Brown, L.D. (1986). {\it Foundations of Exponential Families}. IMS Lecture Notes,
%Monograph Series {\bf 9}, Hayward, California.
%\item \vskip  -2mm DasGupta, A., and Studden, W. (1991). {\it
%Robust Bayesian experimental designs in normal linear models}.
%Annals of Statistics, {\bf 19}, 1244--1256.

\item\vskip -2mm Boraty\'{n}ska, A. (2002). {\it Posterior regret gamma-minimax estimation in a normal model with asymmetric loss function.} Acta Mathematicae, {\bf 29}, 7--13.

\item\vskip -2mm Boraty\'{n}ska, A. (2006). {\it Robust Bayesian prediction with asymmetric loss function in Poisson model of insurance risk}. Acta Universitatis Lodziensis, Folia Oeconomica, {\bf 196}, 123--138. 

\item\vskip-2mm Diaconis, P. and  Ylvisaker, D. (1979). {\it Conjugate priors for exponential families}. Annals of Statistics, {\bf 7}, 269--281.

\item\vskip -2mm Druilhet, P.  and Pommeret D. (2012). {\it Invariant conjugate analysis for exponential families}. Bayesian Analysis, {\bf 7}, 235--248. 
%
%
%\item  \vskip  -2mm Jafari Jozani, M., Marchand, \'E., and
%Parsian, A. (2006). {\it On estimation with weighted balanced-type
%loss function}. Statistics \& Probability Letters, {\bf 76},
%773--780.

\item \vskip -2mm Gelman, A. (2004). {\it Parameterization and Bayesian modelling}. Journal of the American Statistical Association, {\bf 99}, 537--545.

\item\vskip -2mm G\'omez-D\'eniz, E. (2009). {\it Some Bayesian credibility premiums obtained by using posterior regret gamma-minimax methodology.} Bayesian Analysis, {\bf 4}, 223--242. 

\item\vskip-2mm GutiŽrrez-Pe–a, E. (1992). {\it Expected Logarithmic Divergence for Exponential Families}. In Bayesian Statistics 4 (J.M. Bernardo, J.O. Berger, A.P. Dawid y A.F.M. Smith, eds.) Oxford: University Press,  669--674. 

\item\vskip-2mm GutiŽrrez-Pe–a, E. y Smith, A.F.M. (1995). {\it Conjugate Parametrizations for Natural Exponential Families. } Journal of the American Statistical Association, {\bf 90}, 1347--1356. 

\item\vskip  -2mm {Jafari Jozani, M.}, and Parsian, A.
(2008). \textit{Posterior regret $\Gamma$-minimax estimation and
prediction based on $k$-record data under  entropy loss function}.
Communications in Statistics: Theory and Methods, 37, {\bf 14},
2202--2212.

%\item\vskip -2mm Meczarski, M. (1993). {\it Stability  and conditional $\Gamma$-minimaxity in
%Bayesian inference}. Applicationes Mathematicae, {\bf 22}, 117--122.

%\item\vskip -2mm Meczarski, M., and Zielinski, R.  (1991). {\it Stability of Bayesian estimator of
%the Poisson mean under the inexactly specified Gamma prior}. Statistics \& Probability
%Letters, {\bf 12}, 329--333.

\item \vskip  -2mm Robert, C.P. (1996). {\it Intrinsic loss
functions}. Theory  and Decision, {\bf 40}, 192--214.

\item\vskip -2mm  Rios Insua, D., and Ruggeri, F. (2000). {\it Robust Bayesian analysis}.
Lecture Notes in Statistics {\bf 152}, Springer-Verlag, New York.

\item\vskip -2mm  Rios Insua, D., Ruggeri, F., and Vidakovic, B.
(1995). {\it Some results on posterior regret $\Gamma$-minimax
estimation}. Statistics \& Decisions, {\bf 13}, 315--351.

\item \vskip  -2mm Zen, M., and DasGupta, A. (1993). {\it Estimating a binomial parameter: Is robust Bayes
 real Bayes?}  Statistics \& Decisions, {\bf 11}, 37--60.

\end{description}

\end{document}